\documentclass[final,a4paper]{amsart}
\usepackage[latin1]{inputenc}
\usepackage[english]{babel}
\usepackage{amsmath, amsfonts,amssymb,amsopn,amscd,amsthm}
\usepackage{mathpazo}         %
\usepackage{mathrsfs}
\usepackage[final]{graphicx}
\usepackage{indentfirst}
\usepackage{dsfont}
\usepackage{color}
\newtheorem{theorem}{Theorem}[section]
\newtheorem{lemma}[theorem]{Lemma}

\title[Estimating the transition matrix of a Markov chain observed at random times]{Estimating the transition matrix of a Markov chain observed at random times}
\date{\today}

\author{F. Barsotti}\address{FB is with Risk Methodologies, Group Financial Risks, Group Risk Management, UniCredit S.p.A, 20154 Milano. Note that the views presented in this paper are solely those of the author and do not necessarily represent those of UniCredit Spa.}
\email{flavia.barsotti@unicredit.eu}

\author{Y. de Castro}\address{YdC is with Laboratoire de Mathématiques d'Orsay, Universit\'e Paris-Sud, Facult\'e des Sciences d'Orsay, 91405 Orsay, France.}
\email{yohann.decastro@math.u-psud.fr}

\author{T. Espinasse}\address{TE is with Institut Camille Jordan, Université Claude Bernard Lyon 1, 43 boulevard du 11 novembre 1918, 69622 Villeurbanne cedex, France.}
\email{thibault.espinasse@math.univ-lyon1.fr}

\author{P. Rochet}\address{PR is with Laboratoire de Mathématiques Jean Leray, Universit\'e de Nantes, 2 Rue de la Houssinière, 44322 Nantes Cedex 03, France.}
\email{paul.rochet@univ-nantes.fr}

\newcommand{\vect}{\operatorname{vec}}
\newcommand{\id}{\operatorname{I}}
\newcommand{\supp}{\operatorname{supp}}

\newcommand{\lin}{\operatorname{lin}}
\newcommand{\Com}{\operatorname{Com}}
\newcommand{\im}{\operatorname{Im}}

\begin{document}

\maketitle

\normalsize

\begin{abstract}
In this paper we develop a statistical estimation technique to recover the transition kernel $P$ of a Markov chain $X=(X_m)_{m \in \mathbb N}$ in presence of censored data. We consider the situation where only a sub-sequence of $X$ is available and the time gaps between the observations are iid random variables. Under the assumption that neither the time gaps nor their distribution are known, we provide an estimation method which applies when some transitions in the initial Markov chain $X$ are known to be unfeasible. A consistent estimator of $P$ is derived in closed form as a solution of a minimization problem. The asymptotic performance of the estimator is then discussed in theory and through numerical simulations.
\end{abstract}

\maketitle 

\section{Introduction}
\label{sec:int}

Discrete Markov chains are one of the most widely used probabilistic framework for analyzing sequence data in a huge range of application fields. Statistical inference in a Markovian environment has been studied intensively in the literature, giving rise to the definition of various models such as multiple Markov chains \cite{MR0084903,MR0039207}, hidden Markov processes \cite{MR0202264,rabiner1989tutorial}, random walks on graphs \cite{gkantsidis2004random} or renewal processes \cite{smith1958renewal} to cite a few. 

\subsection{Problem}
In this paper we propose a statistical methodology to estimate the transition matrix $P$ from a sequence of censored data. A simple homogenous Markov chain $X = (X_m)_{m \in \mathbb N}$ is observed at random times $T_1,...,T_n$ so that the only available observations consist in a sub-sequence $Y_k:= X_{T_k}$ of the initial process. The time gaps $\tau_k:= T_k - T_{k-1}$ (i.e. the number of jumps) between two consecutive observations are assumed to be positive, independent and identically distributed. \\

\noindent
{\bf Problem}: {\it Can we estimate the transition matrix $P$ of the initial chain $X$ when neither the time gaps $\tau_k$ nor their distribution $\mu$ are known?} \\

\noindent
Without any additional information on the transition kernel $P$, the problem is clearly not identifiable. The novelty of our approach lies in solving this identifiability issue by assuming that some transitions in the initial Markov chain $X$ are known to be unfeasible, that is, the support of $P$ is contained in some maximal set $S$ known to the user. However, even with this new information, the identifiability is only ensured for specific values of $S$. We show for instance that if $S$ contains the entire diagonal or if it is the support of a full bipartite graph, the problem is never identifiable, regardless of the distribution $\mu$. 

A key element in this framework is that the chain $(Y_k)_{k \in \mathbb N}$ remains Markovian, with transition matrix $Q$ that can be expressed as an analytic function of $P$. It results that the problem is identifiable if $P$ is the unique stochastic matrix with compatible support such that there exists an analytic function $f$ verifying $Q = f(P)$. Even when assuming the identifiability, finding a consistent estimation procedure is not straightforward. Using standard non parametric techniques to estimate $f$ seems to be a dead-end, but emphasizes the fact that $f$ and $P$ should be estimated simultaneously. Actually, the main tractable property is that $P$ and $Q$ share the same eigenvectors, although finding the eigenvalues of $P$ from the observations $Y_1, \cdots, Y_N$ remains difficult, as the information on the support of $P$ can not be easily transposed onto conditions on the eigen-elements.

\subsection{Main result}
Our contribution can be described as follows: {\it We estimate the transition matrix $P$ using only the commutativity between $P$ and $Q$. We build an estimator by minimizing the $\ell_{2}$-nom of the Lie bracket with respect to the empirical estimate $\hat Q$ and provide an explicit formula. Moreover, we show the asymptotic normality of the estimator and compute its asymptotic variance. A Monte Carlo simulation study is provided to test its performance, with convincing results.}

To illustrate this model, consider a continuous time Markov chain $Z=(Z_t)_{t>0}$ observed at a discrete time grid $t_1 < ... < t_n$.  In this situation, let $X$ represent the jump process of $Z$ and $\tau_k$ denote the number of jumps occurring between two consecutive observations $Y_k:=Z_{t_k}$. If the discrete time grid is chosen independently from the chain, the time gaps $\tau_k$ are independent random variables {\sl{unknown}} to the practitioner. Their distribution $\mu$ is Poisson in case of a uniform time grid, but one can easily imagine a more involved situation in which $t_1, ...t_n$ are subject to unwanted random effects. The described framework can be found in numerous application fields, since Markov chains are widely recognized for providing faithful representations of real phenomena such as chemical reactions \cite{anderson2011continuous}, financial markets \cite{MR1822778} or waiting lines in queuing theory \cite{gaver1959imbedded}. Markov chains observed at non-regular time intervals are also used for medical studies in \cite{craig2002estimation} to describe the progression of a disease. More general applications of time-varying Markov processes are \cite{MR668189} and \cite{MR0468086}. 

In our setting, the restricted support of $P$ plays a key role in the estimation methodology. Among modern literature contributions, sparsity has become a major interest for statistical inference as it generally provides a significant amount of information that is difficult to fully exploit (see \cite{stewart2009probability,jaaskinen2013sparse}). Here, the sparsity issue is addressed in a specific setting in which the location of some zero entries of $P$ is known. While this considerably simplifies the estimation problem if compared to a framework in which no information is available about the support of $P$, it remains nonetheless a reasonable assumption for most real applications. The spirit of this paper is to present a new approach for inference of sparse Markov transition kernels, as well as to provide a starting point to develop more sophisticated techniques to fully exploit the sparsity in Markov models.

\subsection{Paper organization}
The paper is organized as follows. Section \ref{sec:prob} gives an overview of the statistical framework, describes the estimation problem in detail, and discuss the identifiability issues. Section \ref{sec:est} shows how to characterize and build the estimator of the transition kernel and discusses its asymptotic properties. Section \ref{sec:cs} supports the study with numerical results from a Monte Carlo simulation analysis. Proofs and technical lemmas for our results are gathered in the Appendix.

\section{The problem}
\label{sec:prob}

We consider an irreducible homogenous Markov chain $X=(X_m)_{m \in \mathbb N}$ with finite state space $\mathcal E = \{1,...,N\}$, $N \geq 3$ and transition matrix $P$. We assume that $X$ is observed at random times $T_1,...,T_n$ so that the only available observations consist in the sub-sequence $Y_k:= X_{T_k}, k=1,...,n$. The numbers of jumps $\tau_k:= T_{k-1} - T_k$ between two observations $Y_k$ are assumed to be iid random variables with distribution $\mu$ on $\mathbb N$ and independent from $X$. In this setting, the resulting process $Y=(Y_k)_{k \in \mathbb N}$ remains Markovian in view of the equality:

\begin{align} \mathbb P(Y_{k+1} = j \vert Y_k = i) & = \mathbb P(X_{S_{k+1}} = j \vert X_{S_k} = i) \nonumber \\
& = \sum_{l\geq 0} \mathbb P(X_{S_{k}+l} = j, \tau_{k+1}=l \vert X_{S_k} = i) \nonumber \\
& = \sum_{l \geq 0} \mathbb P(X_{l} = j \vert X_0 = i) \mu(l). \nonumber
\end{align}
Let $G_\mu:[-1,1] \to \mathbb R$ denote the generator function of $\mu$, the transition matrix of $Y$ is thus given by
\begin{equation}\label{Q=fP} Q:= G_\mu(P) = \sum_{l\geq 0} P^l \mu(l).   \end{equation}

We are interested in estimating the original transition matrix $P$ from the available observations $Y_1,...,Y_n$. So far, the problem is not identifiable since neither the time gaps $\tau_k$ nor their distribution $\mu$ are known. Nevertheless, this statistical identifiability issue can be overcome by working with a sparse transition kernel $P$. In this case, we assume that some transitions of the initial process $(X_m)_{m \in \mathbb N}$ are known to be unfeasible, that is, there exists a known set $S \subset \mathcal E^2$ such that 
$$ \supp(P) =  \{ (i,j):   P_{ij} \neq 0 \} \subseteq S.   $$
This additional information restrains the set of possible values of $P$ to
$$ \mathcal A(S) := \{ A \in \mathbb R^{N \times N}: A \mathbf 1 = \mathbf 1, \ \supp(A) \subseteq S \},  $$
which is an affine space of dimension $d-N$, with $d$ the size of $S$. Of course, $P$ is also known to have positive entries, although we choose to overlook this information for now, for simplicity. Assuming that $Q$ is known, we may consider as a solution any stochastic matrix $A \in \mathcal A(S)$ such that $Q=G_\nu(A)$ for some distribution $\nu $ on $\mathbb N$. So, it is possible to recover $P$ exactly from $Q$ if $P$ is the only solution in $\mathcal A(S)$. By slightly relaxing this condition, we say that the problem is \textit{identifiable} if $P$ is the only element in $\mathcal A(S)$ that commutes with $Q$, i.e., if
\begin{equation}\label{idcond} \mathcal A(S) \cap \Com(Q) = \{P \},  \end{equation}
where $\Com(Q)$ denotes the commutant of $Q$. As illustrated in the following lemma, the identifiability of the problem is mainly determined by the value of the support $S$. 

\begin{lemma} \label{identif} The set $\{A \in \mathcal A(S): \  \mathcal A(S) \cap \Com(G_\mu(A)) = \{A \} \}$ is either empty or a dense open subset of $\mathcal A(S)$.
\end{lemma}

This lemma establishes that the problem is either identifiable for almost every possible values of $P$ (with respect to Lebesgue measure) or none, depending on $S$. Remark for instance that the identifiability condition \eqref{idcond} is never verified if $S$ contains the diagonal $\{ (j,j), \ j =1,...,N \}$. Indeed, in this case, the identity matrix $\id$ lies in the intersection $\mathcal A(S) \cap \Com(Q)$ as well as any convex combination $\alpha \id + (1-\alpha) P$ for $\ \alpha \in (0,1)$. Another problematic situation arises if $S$ is the support of a full bipartite graph, resulting in a periodic Markov chain $X$. In this case, the support of $P^3$ is also contained in $S$, which may cause the problem to be non identifiable as soon as $P^3 \neq P$. Similar arguments hold of course for periods other than $2$. Moreover, the problem is not identifiable if $S$ provides insufficient information on $P$. This typically occurs when $d$, the size of $S$, is greater than $ N^2 - N$, or equivalently, if the sparsity degree of $P$ is less than $N$. In this situation, it is easy to show that the dimension of the affine space $\mathcal A(S) \cap \Com(Q)$ is at least $1$, which is obviously incompatible with the identifiability condition given by \eqref{idcond}.\\

While we are able to provide necessary conditions on $S$ for the problem to be identifiable, sufficient conditions turn out to be much harder to obtain. Indeed, this issue involves the companion problem of the eigenvector characterization of weighted directed graphs. Nevertheless, computational study suggests that the combination of the three conditions
\begin{itemize}
	\item $S$ is the support of an aperiodic irreducible Markov chain,
	\item $d \leq N(N - 1)$,
	\item $\exists j, \ (j,j) \notin S$,
\end{itemize}
is sufficient to ensure the almost everywhere identifiability, as we were unable to exhibit a counter-example.\\

To avoid considering critical situations, we will assume throughout the paper that $(X_m)_{m \in \mathbb N}$ is an aperiodic Markov chain. This implies in particular that $P$ has a unique invariant distribution $\pi = (\pi_1,...,\pi_N)$ which is positive for all $i$. Moreover, we assume that the problem is identifiable, i.e., $\mathcal A(S) \cap \Com(Q) = \{P \}$, so that recovering $P$ from the indirect observations $Y_1,...,Y_n$ is achievable.

\section{Construction of the transition kernel estimator}
\label{sec:est}
\noindent We start by introducing some notation. Let $P_0$ be an arbitrary element in $\mathcal A(S)$ and $\phi=(\phi_1 , ... , \phi_{d-N})$ a basis of the difference space
$$ \mathcal A_{\text{lin}}(S) = \mathcal A(S) - \mathcal A(S) = \{ A \in \mathbb R^{N \times N}: A \mathbf 1 = 0, \ \supp(A) \subseteq S \}.  $$
A matrix of the affine space $\mathcal A(S)$ can be decomposed in a unique fashion in function of $P_0$ and $\phi$ as
$$ P_\beta = P_0 + \sum_{j=1}^{d-N} \beta_j \phi_j, $$
for some vector $\beta=(\beta_1,...,\beta_{d-N})^\top \in \mathbb R^{d-N} $. In this setting, the problem of estimating $P$ turns into recovering the corresponding value $\beta$. Consider for convenience its vectorization, which we denote by a small letter, e.g., $p = \vect(P) = (P_{1,1},...,P_{N,1},...,P_{1,N},..., P_{N,N})^\top$. The vector $p$ can be expressed as function of $\beta$ by the relation $ p = p_0 + \Phi \beta$, with $p_0 = \vect(P_0)$ and $\Phi = \left[\vect(\phi_1),...,\vect(\phi_{d-N}) \right]$. When the problem is identifiable, $P$ can be characterized via the Lie bracket with respect to $Q$, as the unique solution in $\mathcal A(S)$ to $\ell(Q,P) = QP - PQ = 0$. Working with the vectorized matrices, the linear operator $p \mapsto \vect[\ell(Q,P)]$ has canonical representation given by $ \Delta(Q) :=  \id \otimes Q - Q^\top \otimes \id$, in view of
$$ \vect( Q P - P Q) = (\id \otimes Q - Q^\top \otimes \id) \vect(P) = \Delta (Q) p. $$
As a result, the information $p = p_0 + \Phi \beta$ and $\Delta(Q) p =  \Delta(Q) [p_0 + \Phi \beta] = 0$ is sufficient to characterize $p$ in this framework. The estimation of $p$ only requires to compute a preliminary estimator of $Q$, say $\hat Q$, which can be directly obtained from the available observations. The most natural choice for $\hat Q$ is arguably the empirical estimator obtained from the state transition frequencies in the sequence $Y_1,...,Y_n$,
\begin{equation}
\displaystyle
\label{eq:Qchap} \hat Q_{ij} = \frac{\sum_{k=1}^{n-1} \mathds 1 \{ Y_k=i,Y_{k+1}=j\}}{ \sum_{k=1}^{n-1} \mathds 1 \{ Y_k=i\} }, \ i,j =1,...,N. \end{equation}
An estimator $\hat p = p_0 + \Phi \hat \beta$ is then quite naturally derived by considering the value $\hat \beta$ for which $\Delta(\hat Q) [p_0 + \Phi \hat \beta]$ is closest to zero. Precisely, we define $\hat \beta$ as a minimizer of
\begin{equation}\label{eqmin}  \beta \mapsto \Vert \Delta(\hat Q) [p_0 + \Phi \beta] \Vert^2 = [p_0 + \Phi \beta]^\top \Delta(\hat Q)^\top \Delta(\hat Q) [p_0 + \Phi \beta].  \end{equation}
If $\Phi^\top \Delta(\hat Q)^\top \Delta(\hat Q) \Phi $ is invertible, the solution is unique, given by
\begin{equation}
\label{eq:Betachap}  \hat \beta = [\Phi^\top \Delta(\hat Q)^\top \Delta(\hat Q) \Phi ]^{-1} \ \Phi^\top \Delta(\hat Q)^\top \Delta(\hat Q)p_0.  \end{equation}
On the contrary, if $\Phi^\top \Delta(\hat Q)^\top \Delta(\hat Q) \Phi $ is singular, we can still define the estimator by picking an arbitrary value among the minimizers. For instance, the solution is obtained via the Moore-Penrose inverse as $  \hat \beta = (\Delta(\hat Q) \Phi)^{\dagger} \Delta(\hat Q)p_0$ (we refer to \cite{MR1408680} for more details on the Moore-Penrose inverse operator). However, nothing indicates that this estimator is close to the true value when $\Delta(\hat Q) \Phi$ is not one-to-one. Actually, the existence of a unique solution is crucial to make the estimator satisfactory. This issue turns out to be closely related to the identifiability of the problem since we can show that condition \eqref{idcond} ensures $\Delta(\hat Q) \Phi$ being of full rank with probability one asymptotically (see Lemma \ref{idcondeq} for a detailed proof). This guarantees that, asymptotically, a unique solution exists. Remark that if the problem is non-identifiable, $\Delta(\hat Q) \Phi$ might be of full rank but not its limit as $n \to \infty$, which would result in a highly unstable, non-consistent estimator. \\

The closed expression of $\hat p := p_0 + \Phi \hat \beta$ enables to derive its asymptotic properties directly from that of $\hat Q$, which we summarize in the next lemma. 

\begin{lemma}\label{hatq} The Markov chain $Y$ is recurrent and share the same invariant distribution $\pi=(\pi_1,...,\pi_N)$ as $X$, which is positive for all $i=1,..,N$. Moreover, $\hat Q$ is unbiased and asymptotically Gaussian with
$$  \forall i,j,k,l = 1,...,N, \   \lim_{n \to \infty} n \operatorname{cov}(\hat Q_{ij}, \hat Q_{kl})  = \left\{ \begin{array}{cl} Q_{ij}(1 - Q_{ij})/\pi_i & \operatorname{if } \ (i,j) = (k,l), \\
 - Q_{ij} Q_{il}/\pi_i & \operatorname{if } \ i=k, \ j\neq l, \\
0 & \operatorname{otherwise.} \end{array} \right. $$
\end{lemma}
\noindent This lemma gathers some well known results on the empirical transition matrix of a finite-state Markov chain. A proof can be found for instance in Theorems 2.7 and 2.15 in \cite{guttorp1995stochastic}. From this result, we deduce that $\hat q = \vect(\hat Q)$ is asymptotically Gaussian, i.e.
\[
\sqrt n (\hat q - q) \overset{d}{\longrightarrow} \mathcal N(0,\Sigma)
\]
for some matrix $\Sigma$ whose expression can be deduced from Lemma \ref{hatq}. We now state our main result.

\begin{theorem}\label{loias} The estimator 
{
\begin{equation}
\label{eq:p} \hat p =  \left[ \id - \Phi [\Phi^\top \Delta(\hat Q)^\top \Delta(\hat Q) \Phi ]^{-1} \ \Phi^\top \Delta(\hat Q)^\top \Delta(\hat Q) \right] p_0 
\end{equation}}
is consistent and asymptotically Gaussian with
$$ \sqrt n (\hat p - p) \overset{d}{\longrightarrow} \mathcal N(0,B \Sigma B^\top),   $$
where $ B = \Phi [ \Phi^\top  \Delta(Q)^\top \Delta(Q) \Phi ]^{-1} \Phi^\top \Delta(Q)^\top \Delta(P)$. 
\end{theorem}

It is worth noting that the value of $\hat p$ does not depend on the initial element $P_0$ nor on the choice of the basis $\phi$. Besides, nothing in the construction of $\hat p$ guarantees that its entries are non-negative. To solve this problem, a natural final step is to consider the stochastic matrix closest to $\hat P$, by vanishing all negative entries and rescaling it so as to obtain an acceptable value. This final solution is clearly a more accurate estimation of $P$. However, we choose to discuss only the properties of the original value $\hat p$ as there are easier to derive and asymptotically equivalent when $S = \supp(P)$.\\

While the proposed transition kernel estimator $\hat{p}$ turns out to be consistent, its efficiency still needs to be discussed. Actually, one can show that $\hat p$ is generally not asymptotically optimal since its limit variance $B \Sigma B^\top$ can be improved. Instead of defining $\hat \beta$ through \eqref{eqmin}, one may consider for instance minimizing a more general quadratic form
\begin{equation}\label{eqmin2}  \beta \mapsto \Vert \Omega \Delta(\hat Q) [p_0 + \Phi \beta] \Vert^2 = [p_0 + \Phi \beta]^\top \Delta(\hat Q)^\top (\Omega^\top \Omega) \Delta(\hat Q) [p_0 + \Phi \beta],  \end{equation}
for some suitably chosen matrix $\Omega \in \mathbb R^{q \times N^2}$, possibly non-square. The only condition we impose on $\Omega$ is that $\Phi^\top  \Delta(Q)^\top (\Omega^\top \Omega) \Delta(Q) \Phi$ must be invertible to guarantee the unicity of the solution, in which case we say that $\Omega$ is admissible. Clearly, the operator $\Omega$ has an influence on the value of the minimizer $\hat \beta_\Omega$, and therefore, on the asymptotic variance of the resulting estimator 
{\begin{equation}
\label{eq:pchap}\hat p_\Omega := p_0 + \Phi \hat \beta_\Omega.
\end{equation}}
By extending the proof of Theorem \ref{loias}, we can show that $\hat p_\Omega$ is asymptotically Gaussian with limit distribution given by

$$ \sqrt n (\hat p_\Omega - p) \overset{d}{\longrightarrow} \mathcal N(0,B(\Omega) \Sigma B(\Omega)^\top),  $$

\noindent for $$B(\Omega) = \Phi [ \Phi^\top  \Delta(Q)^\top (\Omega^\top \Omega) \Delta(Q) \Phi ]^{-1} \Phi^\top \Delta(Q)^\top (\Omega^\top \Omega) \Delta(P).$$
This general approach obviously includes the original procedure corresponding to $\Omega = \id$. Asymptotic optimality can then be derived by aiming for the minimal variance $B(\Omega) \Sigma B(\Omega)^\top$. Using a similar argument as in Proposition 1 in \cite{MR888070}, we show that the minimal variance is reached for any $\Omega$ such that 
$$ (\Omega^\top \Omega) = (\Delta(P) \Sigma \Delta(P)^\top)^\dagger,$$ 
provided that $\Omega$ is admissible (see Lemma \ref{omegaopt} for a detailed proof). This result raises the problem that an optimal value $\Omega^*$ is unknown in practice and has to be estimated beforehand, which can be difficult due to the discontinuity of the Moore-Penrose inversion. Actually, a two-step procedure that consists in plugging-in an estimate $\hat \Omega$ of $\Omega^*$ in \eqref{eqmin2} to compute $\hat \beta_{\hat \Omega}$ might work well in some cases, although theoretical results regarding its performance requires regularity conditions that are hard to verify in practice. For this reason, we suggest to favor the original procedure of Proposition \ref{loias} which provides a consistent estimator $\hat{p}$ by simple means, under no regularity conditions other than the identifiability one given in \eqref{idcond}. Nevertheless, the performances of the two-step estimator $\hat p_{\hat \Omega}$ compared to $\hat{p}$ in various situations are discussed in the next Section via numerical simulations.

\section{Computational study}
\label{sec:cs}
This section is devoted to a Monte Carlo simulation analysis of the proposed methodology. The computational study aims at veryfing the convergence of the estimator $\hat p$ as well as at evaluating the performances of the two-step estimator $\hat p_{\hat \Omega}$  defined in (\ref{eq:pchap}). As discussed in the previous section, the construction of $\hat p_{\hat \Omega}$ involves a preliminary step, namely the estimation of the optimal scaling $(\Omega^{* \top} \Omega^*) = (\Delta(P) \Sigma \Delta(P)^\top)^\dagger$. While $\Delta(P)$ can naturally be estimated from $\Delta(\hat P)$, it remains to build a consistent estimation of $\Sigma$. Actually, this can be made quite easily from observations $Y_1,...,Y_n$. To begin with, the invariant distribution $\pi$ can be estimated by its empirical version
$$ \forall i=1,...,N, \  \hat \pi_i = \frac 1 n \sum_{k=1}^n \mathds 1 \{ Y_k = i \}.$$
It is well known that the resulting estimator $\hat \pi=(\hat \pi_1,...,\hat \pi_N)$ converges to the invariant distribution as soon as the Markov chain is recurrent, which is the case here. We then obtain a consistent estimator $\hat \Sigma$ by replacing $Q$ and $\pi$ by their empirical counterparts in the expression of $\Sigma$, given in Lemma \ref{hatq}. In the following study, the scaling $\hat \Omega = \hat \Sigma^{\frac 1 2} \Delta(P)^\top (\Delta(\hat P) \hat \Sigma \Delta(P)^\top)^\dagger$ is used for the construction of $\hat p_{\hat \Omega}$.\\

The simulations are performed on three examples, corresponding to different values of $P$. The first example deals with an arbitrary sparse transition matrix $P$ for which the support {is} randomly drawn beforehand. The second example investigates an application of our statistical methodology to a queuing model. Finally, the third example considers hollow matrices, for which all entries but the diagonal are non-zero. In each example the transition matrix $P$ is determined beforehand and fixed for the rest of the study. We denote its support with $S:= \supp(P)$. In each case study we consider three different sample sizes $n=200$, $n=1000$ and $n=5000$ and three different distributions for the times gaps $\tau_i$, namely a binomial, Poisson and geometric distribution, the later defined for positive integers only. The whole estimation experiment is repeated $10^4$ times in each setting with the transition matrix $P$ being fixed. Mean squared errors for the two estimators
{\begin{eqnarray}
\label{eq:Rp}
\operatorname{R}(\hat p) &=& \mathbb E \Vert \hat p - p \Vert^2,\\
\label{eq:Rpchap}
\operatorname{R}(\hat p _{\hat \Omega}) &=& \mathbb E \Vert \hat p _{\hat \Omega} - p \Vert^2,
\end{eqnarray}
}
are approximated by the error average over the $10^4$ Monte Carlo repetitions.

\subsection*{Example 1: Random support} This example deals with an arbitrary sparse transition matrix on a state space of size $5$. The support is randomly drawn from independent Bernoulli variables. $P-$entries are drawn from a uniform distribution on $[0,1]$, then rescaled so that $P$ is a transition kernel. The entries are rounded to $2$ decimal digits for ease of readability. We obtain the following transition {matrix} $P$:
\begin{equation}
\label{ex:p1} P = \left[ \begin{array}{ccccc} 
0 & 0.61 & 0 & 0 & 0.39 \\
0.07 & 0 & 0.48 & 0.27 & 0.18 \\
0.53 & 0 & 0.30 & 0 & 0.17 \\
0.18 & 0.20 & 0.27 & 0.35 & 0 \\
0.20 & 0 & 0.69 & 0 & 0.11 
 \end{array} \right].   
 \end{equation}

For the computational study, we proceed as follows. We start by drawing the time gaps $\tau_1,...,\tau_n$ as iid random variables with a given distribution $\mu$ on $\mathbb N$. We consider different values for $n$ and alternative distributions $\mu$ and repeat the experiment in each setting. We let $S_k= \sum_{i=1}^k \tau_i$ for $k=1,...,n$, and simulate a sequence $X_1,X_2,...,X_{S_n}$ of a Markov chain with transition kernel $P$. We keep only the observations to $Y_k =X_{S_k}$ so that we have a sample of size $n$. The process is repeated until all states appear in the sequence $Y_1,...,Y_n$ (in this way, we work conditionally to the event $\hat \pi_i >0$). From these observations, we build $\hat p$ following its closed expression in Theorem \ref{loias} and the two-step estimator $\hat p_{\hat \Omega}$ defined in \eqref{eq:pchap} using the procedure detailed above. The whole experiment is repeated $10^4$ times with the same value of $P$ for the three different sample sizes $n=200$, $n=1000$ and $n=5000$ and three different distributions, namely a binomial $\mu \sim \mathcal B(5,0.3)$, standard Poisson $\mu \sim \mathcal P(1)$ and geometric distribution $\mu \sim  \mathcal G(0.5)$. Mean squared errors given in (\ref{eq:Rp}) and (\ref{eq:Rpchap}) are evaluated and reported in Table 1 below with standard deviations in brackets. \\

{\scriptsize
\begin{table}[ht]
\begin{center}
\begin{tabular}{|c||c|c|c|c|c|c|c|c|c|} \hline 
$ n $ & \multicolumn{3}{c|}{$200$} &  \multicolumn{3}{c|}{ $1000$} & \multicolumn{3}{c|}{$5000$} \\
\hline
$\mu$ & $\mathcal B(5,0.3)$ & $ \mathcal P(1) $ & $ \mathcal G(0.5)$ & $\mathcal B(5,0.3)$ & $ \mathcal P(1) $ & $\mathcal G(0.5)$ & $\mathcal B(5,0.3)$ & $\mathcal P(1) $ & $\mathcal G(0.5)$ \\
\hline
\hline
$\operatorname{R}(\hat p)$ & $\!\!\! \begin{array}{c} 0.5469 \\ (0.0033) \end{array}\!\!\! $ & $ \!\!\! \begin{array}{c} 0.5189 \\ (0.0030) \end{array}\!\!\!  $ & $\!\!\! \begin{array}{c} 0.3514 \\ (0.0021) \end{array}\!\!\! $ & $\!\!\! \begin{array}{c} 0.1637 \\ (0.0010) \end{array}\!\!\! $ & $ \!\!\! \begin{array}{c} 0.1371 \\ (0.0008) \end{array}\!\!\! $ & $ \!\!\! \begin{array}{c} 0.0835 \\ (0.0004) \end{array}\!\!\! $ & $ \!\!\! \begin{array}{c} 0.0362 \\ (0.0002) \end{array} \!\!\! $ & $\!\!\! \begin{array}{c} 0.0286 \\ (0.0002) \end{array}\!\!\! $ & $\!\!\! \begin{array}{c} 0.0170 \\ (0.0001) \end{array}\!\!\!$     \\
\hline
$\operatorname{R}(\hat p_{\hat \Omega})$ & $\!\!\! \begin{array}{c} 1.2113 \\ (0.0385) \end{array}\!\!\! $ & $ \!\!\! \begin{array}{c} 1.0901 \\ (0.0232) \end{array} \!\!\! $ & $\!\!\! \begin{array}{c} 0.4928 \\ (0.0074) \end{array}\!\!\! $ & $\!\!\! \begin{array}{c} 0.1668 \\ (0.0015) \end{array}\!\!\! $ & $ \!\!\! \begin{array}{c} 0.1389 \\ (0.0011) \end{array}\!\!\! $ & $ \!\!\! \begin{array}{c} 0.0782 \\ (0.0004) \end{array}\!\!\! $ & $\!\!\! \begin{array}{c} 0.0301 \\ (0.0002) \end{array} \!\!\! $ & $\!\!\! \begin{array}{c} 0.0249 \\ (0.0001) \end{array}\!\!\! $ & $\!\!\! \begin{array}{c} 0.0148 \\ (0.0001) \end{array}\!\!\!$ \\
\hline
\end{tabular}
\vspace{0.3cm}
\caption{\footnotesize \textbf{Monte Carlo Experiment Results.} The table contains summary statistics of Monte Carlo simulation results based on $P$ matrix given in Eq. (\ref{ex:p1}).  Three different sample sizes $n$ and three different distributions $\mu$ are considered. Mean squared errors $\operatorname{R}(\hat p)$ and $\operatorname{R}(\hat p _{\hat \Omega})$ defined in Eq. (\ref{eq:Rp})-(\ref{eq:Rpchap}) are reported with the corresponding standard deviations in brackets. Monte Carlo errors are based on $10^4$ repetitions.}
\end{center}
\label{tab:ex1}
\end{table}
}

\normalsize

Theoretical results described in previous sections are now confirmed by the Monte Carlo simulation. For a small sample size ($n=200$), the estimation of $P$ is obviously difficult and it shows a mean squared error $R(\hat{p})=0.35$ in the most favorable case, corresponding to an average squared error of approximately $0.022$ per entry. The two-step procedure considerably deteriorates the estimation for $n=200$, regardless of the distribution of the time gaps. Interesting insights arise for the sample size $n=1000$. In this case, the two estimators $\hat p $ and $\hat p_{\hat \Omega}$ show a comparable performance, with $\hat p$ being slightly better for the binomial and Poisson scenarios, while $\hat p_{\hat \Omega}$ appears to be preferable for the geometric distribution. The transition matrix $P$ is relatively well estimated in this case, especially for geometric times, with an average squared error of approximately $0.005$ per entry. Finally, for a large sample size $n=5000$, the transition matrix $P$ is very well estimated by both methods, with significantly better results for the two-step estimator $\hat p_{\hat \Omega}$. \\

While the distribution $\mu$ seems to have a non negligible impact on the efficiency of the estimation, it is difficult to establish the nature of its influence. The geometric distribution reveals to be the most favorable case here, which was to be expected since it is the only one for which the event $\tau_i =0$ is ruled out. This means that two consecutive observations in the process $Y$ are always different, which is obviously desirable. For the binomial an Poisson cases, it is not rare that the process remains unchanged for two or more consecutive observations of $Y$ which somewhat reduces the number of observations. This explains the better results obtained for the case $\tau \sim \mathcal G(0.5)$ if compared to the other two settings.

\subsection*{Example 2: Queuing model} This example considers the application of our statistical methodology to a queuing model. We want to evaluate the influence of the number of persons in a waiting line at time $t$ on the sub-sequent state of the queue, i.e. at time $t+1$. States represent the number of persons in the queue. For simplicity we assume that the only possible transitions are the arrival or departure of someone. The state of the waiting line is measured only at particular times (e.g. every hour) and the number of transitions $\tau_k$ between two consecutive observations $Y_{k-1}$ and $Y_{k}$ are assumed iid with unknown distribution $\mu$. For convenience, we assume a maximum number of persons in the queue equal to $10$. Thus, the Markov chain $X$ has $11$ possible states and a transition matrix $P$ whose only non-zero entries are $P_{i+1, i}$ and $P_{i, i+1}$ for $i=1,...,10$. These entries are not chosen too far from $0.5$ so that $\hat \pi$ is positive with relatively high probability, even for small sample sizes. The actual transition matrix used for the simulations is the following matrix $P$:

\begin{equation}
\label{ex:p2} P = \left[ \begin{array}{ccccccccccc} 
0 & 1 & 0 & 0 & 0 & 0 & 0 & 0 & 0 & 0 & 0 \\
0.53 & 0 & 0.47 & 0 & 0 & 0 & 0 & 0 & 0 & 0 & 0 \\
0 & 0.65 & 0 & 0.35 & 0 & 0 & 0 & 0 & 0 & 0 & 0 \\
0 & 0 & 0.45 & 0 & 0.55 & 0 & 0 & 0 & 0 & 0 & 0 \\
0 & 0 & 0 & 0.30 & 0 & 0.70 & 0 & 0 & 0 & 0 & 0 \\
0 & 0 & 0 & 0 & 0.62 & 0 & 0.38 & 0 & 0 & 0 & 0 \\
0 & 0 & 0 & 0 & 0 & 0.68 & 0 & 0.32 & 0 & 0 & 0 \\
0 & 0 & 0 & 0 & 0 & 0 & 0.64 & 0 & 0.36 & 0 & 0 \\
0 & 0 & 0 & 0 & 0 & 0 & 0 & 0.52 & 0 & 0.48 & 0 \\
0 & 0 & 0 & 0 & 0 & 0 & 0 & 0 & 0.61 & 0 & 0.39 \\
0 & 0 & 0 & 0 & 0 & 0 & 0 & 0 & 0 & 1 & 0 \\
 \end{array} \right].   
 \end{equation}

Remark that the first and last rows of $P$ are known since they contain only one non-zero element. As in the previous example, we consider three sample sizes $n=200$, $n=1000$ and $n=5000$ and three distributions $\mu = \mathcal B(2,0.5)$, $\mu = \mathcal P(1)$ and $\mu = \mathcal G(0.5)$. Results are gathered in Table 2.

{\scriptsize
\begin{table}[ht]
\begin{center}
\begin{tabular}{|c||c|c|c|c|c|c|c|c|c|} \hline 
$ n $ & \multicolumn{3}{c|}{$200$} &  \multicolumn{3}{c|}{ $1000$} & \multicolumn{3}{c|}{$5000$} \\
\hline
$\mu$ & $\mathcal B(2,0.5)$ & $ \mathcal P(1) $ & $ \mathcal G(0.5)$ & $\mathcal B(2,0.5)$ & $ \mathcal P(1) $ & $\mathcal G(0.5)$ & $\mathcal B(2,0.5)$ & $\mathcal P(1) $ & $\mathcal G(0.5)$ \\
\hline
\hline
$\operatorname{R}(\hat p)$ & $\!\!\! \begin{array}{c} 0.5449 \\ (0.0028) \end{array}\!\!\! $ & $ \!\!\! \begin{array}{c} 0.6530 \\ (0.0035) \end{array}\!\!\!  $ & $\!\!\! \begin{array}{c} 0.4527 \\ (0.0026) \end{array}\!\!\! $ & $\!\!\! \begin{array}{c} 0.1763 \\ (0.0013) \end{array}\!\!\! $ & $ \!\!\! \begin{array}{c} 0.2296 \\ (0.0016) \end{array}\!\!\! $ & $ \!\!\! \begin{array}{c} 0.1253 \\ (0.0009) \end{array}\!\!\! $ & $ \!\!\! \begin{array}{c} 0.0421 \\ (0.0004) \end{array} \!\!\! $ & $\!\!\! \begin{array}{c} 0.0548 \\ (0.0005) \end{array}\!\!\! $ & $\!\!\! \begin{array}{c} 0.0258 \\ (0.0002) \end{array}\!\!\!$     \\
\hline
$\operatorname{R}(\hat p_{\hat \Omega})$ & $\!\!\! \begin{array}{c} 1.0813 \\ (0.0059) \end{array}\!\!\! $ & $ \!\!\! \begin{array}{c} 1.2287 \\ (0.0064) \end{array} \!\!\! $ & $\!\!\! \begin{array}{c} 0.8967 \\ (0.0043) \end{array}\!\!\! $ & $\!\!\! \begin{array}{c} 0.3176 \\ (0.0036) \end{array}\!\!\! $ & $ \!\!\! \begin{array}{c} 0.3788 \\ (0.0038) \end{array}\!\!\! $ & $ \!\!\! \begin{array}{c} 0.2657 \\ (0.0024) \end{array}\!\!\! $ & $\!\!\! \begin{array}{c} 0.0211 \\ (0.0002) \end{array} \!\!\! $ & $\!\!\! \begin{array}{c} 0.0440 \\ (0.0004) \end{array}\!\!\! $ & $\!\!\! \begin{array}{c} 0.0284 \\ (0.0003) \end{array}\!\!\!$ \\
\hline
\end{tabular}
\vspace{0.3cm}
\caption{\footnotesize\textbf{Monte Carlo Experiment Results.} The table contains summary statistics of Monte Carlo simulation results based on $P$ matrix given in Eq. (\ref{ex:p2}).  Three different sample sizes $n$ and three different distributions $\mu$ are considered. Mean squared errors $\operatorname{R}(\hat p)$ and $\operatorname{R}(\hat p _{\hat \Omega})$ defined in Eq. (\ref{eq:Rp})-(\ref{eq:Rpchap}) are reported with the corresponding standard deviations in brackets. Monte Carlo errors are based on $10^4$ repetitions.}
\end{center}
\label{tab:exemple2}
\end{table}
}

\normalsize
Although the number of states is more than doubled compared to the previous example, the Monte Carlo simulation show similar results. This is due to the fact that the difficulty in estimating $P$ is mostly determined by its number of non trivial entries, rather than by its dimension. These values are quite similar in both example ($18$ non trivial entries in this example against $16$ in the previous one). The first assumption for $\mu$, namely the Binomial distribution $\mathcal B(2,0.5)$, leads to a sparse empirical transition matrix $Q$. Indeed, in this case, $Q$ is a convex combination of $\id$, $P$ and $P^2$ which implies that its non-zero entries are at a distance of at most $2$ from the main diagonal. As a result, the estimation $\hat Q$ is somehow more accurate in this case compared to a situation in which all state transitions are possible in the chain $Y$. On the other hand, the high probability of observing the same realization at two consecutive times (due to $\mu(0) =  0.25$) deteriorates the estimation of $P$. Nevertheless, the most favorable case remains the geometric distribution for all sample sizes.\\

Similar conclusions can be drawn regarding the relative efficiency of $\hat p_{\hat \Omega}$ and $\hat p$, as it appears clearly that $\hat p_{\hat \Omega}$ outperforms $\hat p$ only when a large number of observations are available. Interestingly, $\hat p$ remains significantly better even for $n=5000$ in the geometric scenario.

\subsection*{Example 3: Hollow matrix} This final example deals with a transition matrix $P$ with zero diagonal, sometimes referred to as hollow matrix. This case corresponds to a Markov chain $X$ that necessarily changes state at each transition. The matrix $P$ used for the simulation is the following:
\begin{equation}
\label{ex:p3} P = \left[ \begin{array}{cccc} 
0 & 0.22 & 0.33 & 0.45 \\
0.38 & 0 & 0.06 & 0.56 \\
0.40 & 0.13 & 0 & 0.47 \\
0.42 & 0.20 & 0.38 & 0 
 \end{array} \right].   
 \end{equation}
Following the same simulation structure described in the two previous examples, we consider alternative sample sizes and distributions of time gaps. Results are summarized in Table 3.

{\scriptsize
\center
\begin{table}[http!]
\begin{center}
\begin{tabular}{|c||c|c|c|c|c|c|c|c|c|} \hline 
$ n $ & \multicolumn{3}{c|}{$200$} &  \multicolumn{3}{c|}{ $1000$} & \multicolumn{3}{c|}{$5000$} \\
\hline
$\mu$ & $\mathcal B(2,0.5)$ & $ \mathcal P(1) $ & $ \mathcal G(0.5)$ & $\mathcal B(2,0.5)$ & $ \mathcal P(1) $ & $\mathcal G(0.5)$ & $\mathcal B(2,0.5)$ & $\mathcal P(1) $ & $\mathcal G(0.5)$ \\
\hline
\hline
$\operatorname{R}(\hat p)$ & $\!\!\! \begin{array}{c} 0.3386 \\ (0.0036) \end{array}\!\!\! $ & $ \!\!\! \begin{array}{c} 0.3838 \\ (0.0037) \end{array}\!\!\!  $ & $\!\!\! \begin{array}{c} 0.2879 \\ (0.0032) \end{array}\!\!\! $ & $\!\!\! \begin{array}{c} 0.1501 \\ (0.0025) \end{array}\!\!\! $ & $ \!\!\! \begin{array}{c} 0.1978 \\ (0.0029) \end{array}\!\!\! $ & $ \!\!\! \begin{array}{c} 0.1093 \\ (0.0020) \end{array}\!\!\! $ & $ \!\!\! \begin{array}{c} 0.0524 \\ (0.0016) \end{array} \!\!\! $ & $\!\!\! \begin{array}{c} 0.0799 \\ (0.0021) \end{array}\!\!\! $ & $\!\!\! \begin{array}{c} 0.0271 \\ (0.0008) \end{array}\!\!\!$     \\
\hline
$\operatorname{R}(\hat p_{\hat \Omega})$ & $\!\!\! \begin{array}{c} 0.3561 \\ (0.0039) \end{array}\!\!\! $ & $ \!\!\! \begin{array}{c}  0.4011 \\ (0.0039) \end{array} \!\!\! $ & $\!\!\! \begin{array}{c} 0.3003 \\ (0.0035) \end{array}\!\!\! $ & $\!\!\! \begin{array}{c} 0.1646 \\ (0.0028) \end{array}\!\!\! $ & $ \!\!\! \begin{array}{c} 0.2170 \\ (0.0032) \end{array}\!\!\! $ & $ \!\!\! \begin{array}{c} 0.1167 \\ (0.0022) \end{array}\!\!\! $ & $\!\!\! \begin{array}{c} 0.0590 \\ (0.0019) \end{array} \!\!\! $ & $\!\!\! \begin{array}{c} 0.0880 \\ (0.0024) \end{array}\!\!\! $ & $\!\!\! \begin{array}{c} 0.0307\\ (0.0011) \end{array}\!\!\!$ \\
\hline
\end{tabular}
\vspace{0.3cm}
\end{center}
\caption{ \footnotesize \textbf{Monte Carlo Experiment Results.} The table contains summary statistics of Monte Carlo simulation results based on $P$ matrix given in Eq. (\ref{ex:p3}).  Three different sample sizes $n$ and three different distributions $\mu$ are considered. Mean squared errors $\operatorname{R}(\hat p)$ and $\operatorname{R}(\hat p _{\hat \Omega})$ defined in Eq. (\ref{eq:Rp})-(\ref{eq:Rpchap}) are reported with the corresponding standard deviations in brackets. Monte Carlo errors are based on $10^4$ repetitions.}
\label{tab:ex3}
\end{table}
}

\normalsize
In this example, $\hat p$ and $\hat p_{\hat \Omega}$ show comparable performances for all sample size, with yet slightly better results for $\hat p$. Surprisingly, the theoretical asymptotic results seem to not be verified even for sample as large as $n=5000$. A dedicated simulation for a sample size $n=10000$ has been performed and shows that the mean squared error of $\hat p_{\hat \Omega}$ does eventually become smaller than that one  of $\hat p$, however this occurs for very large $n$, in some sense confirming that the regular procedure should be favored in most practical situations.\\

The Monte Carlo experiment is performed on three different situations allowing to draw similar conclusions. We have observed the convergence of both estimators $\hat p$ and $\hat p_{\hat \Omega}$ to  the true value $p$ in all considered examples. Moreover, the two-step procedure to construct the asymptotically estimator $\hat p_{\hat \Omega}$ has appeared unsatisfactory in most cases, with a significant improvement with respect $\hat p$ only for some cases with large samples (from $n=1000$ or even $n>5000$ in the last example). These simulations confirm the theoretical results as well as the conclusion that the regular estimator $\hat p$ must be preferred for both its stability and easiness of implementation. However, while the performances of $\hat p_{\hat \Omega}$ are disappointing, we observe that the scaling $\hat \Omega$ used for its construction can be considered as a naive estimation of the theoretical optimal scaling. The small sample properties of $\hat p_{\hat \Omega}$ can perhaps be improved by using different estimation techniques for $\Omega^*$, although this has not been investigated. 

\section{Conclusion}
\label{sec:conc}
This paper investigates the problem of estimating the transition kernel $P$ of a discrete Markov chain in presence of censored data, i.e. when only a sub-sequence of the chain is observable. The original contribution is the development and proposal of a statistical methodology to recover the transition matrix $P$ when the time intervals between two observations are random, iid and unknown. To overcome the identifiability issue in this setting, $P$ is assumed to be sparse with its zeroes location partially known.
The novelty of our approach lies in the role played by the sparsity of $P$ when working in Markovian models with censored data, since the proposed methodology is able to capture and exploit the information content behind this setting. In the framework studied in this paper, the available observations consist of a Markov chain with a transition matrix $Q$ that commutes with $P$. Once characterized the transition matrix $P$ via the Lie bracket with respect to $Q$, we show how to build an estimator by means of the empirical transition matrix of the observations. A consistent estimator $\hat{p}$ is given in closed form as function of $\hat{Q}$ and its asymptotic properties are derived. Theoretical results are supported by a Monte Carlo simulation study to verify the convergence of the estimator and analyze its performance in various situations. 

In this paper, we focus on a situation where the support of the transition matrix $P$ is partially known. While this assumption turns out to be important to make the problem feasible, we can imagine practical cases in which $P$ is sparse with unknown support. Numerous questions arise in this context, such as recovering the minimal support for a matrix in the commutant of $Q$ or determining necessary and sufficient conditions for the problem to be identifiable. We are optimistic that the current paper provides a significant starting point to tackle these questions in future research works.

\section{Appendix}\label{sec:app}

\subsection{Technical lemmas}

\begin{lemma}\label{idcondeq} The problem is identifiable if, and only if, $\Delta(Q) \Phi$ is of full rank.
\end{lemma}

\noindent \textit{Proof.} Writing the sets $\mathcal A(S) $ and $\Com(Q)$ as the Minkowsky sums $\mathcal A(S) = \{ P \} + \mathcal A_{\lin}(S)$ and $\Com(Q) = \{P \} + \Com(Q)$, we deduce
$$\mathcal A(S) \cap \Com(Q) = \{ P \} + \mathcal A_{\lin}(S) \cap \Com(Q).$$ 
Thus, the identifiability condition is equivalent to $\mathcal A_{\lin}(S) \cap \Com(Q) = \{ 0 \}$. Since $\Phi$ is of full rank, $\ker[\Delta(Q) \Phi] = \{ 0 \}$ holds if, and only if, $\ker[\Delta(Q)] \cap \im(\Phi) = \{ 0 \}$, where $\im$ denotes the image. The result follows by pointing out that $\ker[\Delta(Q)] = \vect[\Com(Q)] : = \{ \vect(A) : A \in \Com(Q) \}$ and $\im(\Phi) = \vect[\mathcal A_{\lin}(S)]$.\\

\begin{lemma}\label{omegaopt} If $ \Phi^\top \Delta(Q)^\top \left(\Delta(P) \Sigma \Delta(P)^\top \right)^\dagger \Delta(Q) \Phi$ is invertible, any matrix $\Omega^*$ such that 
$$ \Omega^{*\top} \Omega^* = \left(\Delta(P) \Sigma \Delta(P)^\top \right)^\dagger,   $$
is asymptotically optimal in the sense that $ B(\Omega) \Sigma B(\Omega)^\top - B(\Omega^*) \Sigma B(\Omega^*)^\top$ is positive semi-definite for all admissible $\Omega$.
\end{lemma}

\noindent \textit{Proof.} Let $\Omega^*$ be such a matrix and let $D= \Omega^* \Delta(Q) \Phi$. The operator $\id - D (D^\top D)^{-1} D^\top$ is an orthogonal projector and is therefore positive semi-definite. Let $\Omega$ be admissible and 
$$C = \left[ \Phi^\top \Delta(Q)^\top (\Omega^\top \Omega) \Delta(Q) \Phi \right]^{-1} \Phi^\top \Delta(Q)^\top \Delta(P) \Omega^{* \top}.$$ 
We know that $ C C^\top - C D (D^\top D)^{-1} D^\top C^\top$ is also positive semi-definite, which yields the wanted result.

\subsection{Proofs} 
$$  $$

\vspace{-0.3cm}

\noindent \textbf{Proof of Lemma \ref{identif}.} Recall that $\mathcal A(S) = \{ P_\beta = P_0 + \sum_{j=1}^{d-N} \beta_j \phi_j: \beta \in \mathbb R^{d-N} \}$. From Lemma \ref{idcondeq}, we know the problem is identifiable if, and only if, $\Delta(G_\mu(P)) \Phi$ is of full rank, i.e. if 
$$  \operatorname{det}(\Phi^\top \Delta(G_\mu(P))^\top \Delta(G_\mu(P)) \Phi) \neq 0.$$ 
Since the map $g:\beta \mapsto \operatorname{det}(\Phi^\top \Delta(G_\mu(P_\beta))^\top \Delta(G_\mu(P_\beta)) \Phi)$ is analytic, $g^{-1}(\{0\})$ is either equal to $\mathbb R^{d-N}$ or is a nowhere dense closed subset of $\mathbb R^{d-N}$.\\

\noindent \textbf{Proof of Theorem \ref{loias}.} If the problem is identifiable, then $\text{ker}[\Delta(Q) \Phi] = \{ 0 \}$ and the map
$$ F: A \mapsto  \left[ I - \Phi \left[ \Phi^\top \Delta(A)^\top \Delta(A) \Phi \right]^{-1} \Phi^\top \Delta(A)^\top \Delta(A) \right] p_0 $$
is continuously differentiable at $A = Q$. Since $\hat Q$ converges in probability to $Q$, we get by Cramer's theorem
\begin{equation}\label{delta} \sqrt n (\hat p - p) = \sqrt n (F(\hat Q) - F(Q)) = \sqrt n \ \nabla F_Q(\hat Q - Q) + o_P(1), \end{equation} 
where $\nabla F_Q$ denotes the differential of $ F$ at $Q$. Direct calculation gives for $H \in \mathbb R^{N \times N}$,
\begin{align} \nabla F_Q(H) = & \lim_{t \to 0} \frac{F(Q + t H) - F(Q)}t \nonumber \\
 = & - \Phi \left[ \Phi^\top  \Delta(Q)^\top \Delta(Q) \Phi \right]^{-1} \Phi^\top  \left[\Delta(Q)^\top \Delta(H) +  \Delta(H)^\top \Delta(Q) \right] p. \nonumber
\end{align}
Noticing that $\Delta(Q) p = 0$ and $\Delta(H) p = - \Delta(P) h$, we get
$$  \nabla F_Q(H) = \Phi \left[ \Phi^\top  \Delta(Q)^\top \Delta(Q) \Phi \right]^{-1} \Phi^\top \Delta(Q)^\top \Delta(P) h.   $$
We now use that $\sqrt n (\hat q - q)  \overset{d}{\longrightarrow} \mathcal N(0, \Sigma)$ combined with \eqref{delta} to complete the proof.

\bibliographystyle{acm}
\bibliography{refs}

\end{document}